\documentclass[11pt,doublespaced]{amsart}
\usepackage{latexsym,amsfonts,amsmath,amssymb}
\usepackage[top=1.4in, bottom=1.4in, left=1.4in, right=1.4in]{geometry}

\newcommand{\IR}{{\mathbb{R}}}
\newcommand{\ID}{{\mathbb{D}}}
\newcommand{\ZZ}{{\mathbb{Z}}}
\newcommand{\IC}{{\mathbb{C}}}

\newcommand{\al}{{\alpha}}

\newcounter{smalllist}

\newtheorem{theorem}{Theorem}

\newtheorem{lemma}{Lemma}[section]
\newtheorem{prop}[lemma]{Proposition}

\theoremstyle{definition}

\newtheorem{remark}[lemma]{Remark}

\let\llldots=\ldots
\def\ldots{\llldots{}}

\numberwithin{equation}{section}

\begin{document}

\title[A note on Poisson brackets for OPUC]{A note on Poisson brackets for orthogonal polynomials on the unit circle}
\author{Irina Nenciu}
\address{Irina Nenciu\\
         Department of Mathematics, Statistics and Computer Science\\
         University of Illinois at Chicago\\
         Chicago, IL 60607, USA \textit{and} Institute of Mathematics ``Simion Stoilow'' of the Romanian Academy\\ Bucharest\\ Romania}
\email{nenciu@uic.edu}
\thanks{The author wishes to thank Maria Cantero, Barry Simon, and Takayuki Tsuchida for their useful comments and suggestions.
This research was supported by NSF grants DMS-0111298 and DMS-0701026. Part of this work was
done while the author was a member of the School of Mathematics of
the Institute for Advanced Study in Princeton.}

\begin{abstract}
The connection of orthogonal polynomials on the unit circle (OPUC) 
to the defocusing Ablowitz-Ladik integrable system
involves the definition of a Poisson structure on the space of Verblunsky coefficients.
In this paper, we compute the complete set of Poisson brackets for the monic orthogonal and the orthonormal
polynomials on the unit circle, as well as for the second kind polynomials and the Wall polynomials,
This answers a question posed by Cantero and Simon, ~\cite{CanSim}, for the case of measures
with finite support. We also show that the results hold for the case of measures with
periodic Verblunsky coefficients.
\end{abstract}

\maketitle

\section{Introduction and background}

Over the past few years, there has been a flurry of activity concerning the defocusing Ablowitz-Ladik system and
its connection to the theory of orthogonal polynomials on the unit circle (OPUC). Among other things, this involves the definition of a Poisson structure on the set of Verblunsky coefficients.
Let us elaborate:
Given a probability measure $d\mu$ on $S^1$, the
unit circle in $\IC$, we can construct an orthonormal system of
polynomials, $\phi_n$, by applying the Gram--Schmidt procedure to
$\{1,z,\ldots\}$.  These obey a recurrence relation; however, to
simplify the formulae, we will present the relation for the monic
orthogonal polynomials $\Phi_n$:
\begin{equation}\label{E:PhiRec}
\Phi_{n+1}(z) = z\Phi_n(z)   - \bar\alpha_n \Phi_n^*(z),
\end{equation}
\begin{equation}\label{E:Phi*Rec}
\Phi_{n+1}^*(z)=\Phi_n^*(z)-\al_nz\Phi_n(z).
\end{equation}
Here, $\alpha_n\in\IC$ are recurrence coefficients, which are called
Verblunsky coefficients, and $\Phi_n^*$ denotes the reversed
polynomial, $\Phi_n^*(z)=z^n\overline{\Phi_n\left(\tfrac{1}{\bar
z}\right)}$. Notice that $\al_n=-\overline{\Phi_{n+1}(0)}.$
When $d\mu$ is supported at exactly $N$ points,
$\alpha_n\in\ID=\{z\in\IC\,|\,|z|<1\}$ for $0\leq n\leq {N-2}$ while
$\alpha_{N-1}$ is a unimodular complex number. If the support of
$d\mu$ is infinite, then the recurrence formula \eqref{E:PhiRec}
will produce an infinite sequence $\{\al_n\}_{n\geq0}$ of Verblunsky
coefficients, all of whom are inside the unit disc $\ID$. In both of
these cases, there is a 1-1 correspondence between the measure and
the sequence of coefficients. For more details on the theory
of polynomials orthogonal on the unit circle, we refer the reader to
the treatise of B. Simon, \cite{Simon1} and \cite{Simon2}.

In the finite case ($d\mu$ supported at $N$ points), one can introduce a Poisson bracket on the space of Verblunsky coefficients,
which, for two smooth functions $f$ and $g$ on $\ID^{N-1}\times S^1$, is\footnote{Note that the sum defining this Poisson bracket
goes only up to $N-2$ and hence the last Verblunsky coefficient, $\al_{N-1}\in S^1$, is a Casimir.}
\begin{equation}\label{E:PbDefn}
\{f,g\}=i\sum_{k=0}^{N-2} \rho_k^2 \left[\frac{\partial f}{\partial \bar\al_k}\frac{\partial g}{\partial \al_k}
-\frac{\partial f}{\partial \al_k}\frac{\partial g}{\partial \bar\al_k}\right],
\end{equation}
where $\rho_k=\sqrt{1-|\al_k|^2}$ and, for a complex variable $\al=u+iv$, $u,v\in\IR$, the partial derivatives are defined
as usual by
\begin{equation*}
\frac{\partial}{\partial\al}=\frac12\left[\frac{\partial}{\partial u}-i\frac{\partial}{\partial v}\right],
\qquad
\frac{\partial}{\partial\bar\al}=\frac12\left[\frac{\partial}{\partial u}+i\frac{\partial}{\partial v}\right].
\end{equation*}
If the sequence of Verblunsky coefficients is infinite, it is very challenging in general to define
a Poisson bracket. Nonetheless, we can also restrict ourselves to the case when the Verblunsky coefficients are periodic,
with period $p$. We assume this period to be even, and we recall (see \cite{NenSim} and \cite{Nen2}) that in this situation
there exists a Poisson bracket on the (finite dimensional!) space
$$
\ID^p=\bigl\{(\alpha_0,\dots,\alpha_{p-1})\,|\,\alpha_j\in\ID\bigr\}\,.
$$
This bracket is given by
\begin{equation}\label{E:PbDefnPeriodic}
\{f,g\}=i\sum_{k=0}^{p-1} \rho_k^2 \left[\frac{\partial f}{\partial \bar\al_k}\frac{\partial g}{\partial \al_k}
-\frac{\partial f}{\partial \al_k}\frac{\partial g}{\partial \bar\al_k}\right]\,,
\end{equation}
and, unlike in the previous case, it is also non-degenrate (i.e. it has no Casimirs).

\begin{remark}
Let us note here that, while using the same symbol for the two Poisson brackets in \eqref{E:PbDefn} and \eqref{E:PbDefnPeriodic}
is an abuse of notation, it will not lead to any confusion in our results here. 
Indeed, the purpose of this note is to compute the Poisson brackets of the
monic orthogonal and orthonormal polynomials, thought of as
functions of the Verblunsky coefficients, as well as of other
polynomials naturally associated to the problem, e.g. second kind
polynomials and Wall polynomials (for definitions, see the second part of the Introduction). One of the main properties of
these polynomials, as functions of the $\alpha$'s, is that each polynomial depends only on a finite 
number of Verblunsky coefficients, $\alpha_0$ through $\alpha_k$, with the largest index $k$ being
uniquely determined by the degree of the polynomial. In particular, this means that as long as this largest
index $k$ is strictly smaller than either $N-1$ or $p$, the Poisson brackets for these polynomials
are given by identical expressions in the finite and the periodic cases. 
\end{remark}

The first to compute Poisson brackets
in the context presented here were, to the best of our knowledge, Nenciu and Simon, \cite{NenSim}.
Their purpose was to understand the complete integrability of the
defocusing Ablowitz-Ladik (AL) equation with periodic boundary conditions;
in the process, they compute certain fairly complicated
combinations of the Poisson brackets of the Wall polynomials--see Proposition~\ref{P:cWP}.
(For the origin of the Ablowitz-Ladik equation, see \cite{AL1}, \cite{AL2}; for recent results
on this integrable system, obtained through its connection to OPUC, see for example
\cite{CanSim, GekNen, KilNen, LiNen, Nen}, or the review paper
\cite{Nen3}.) In \cite{CanSim}, Cantero and Simon found all
but one of the Poisson brackets for the monic orthogonal and the second kind
polynomials. In this paper, we give all the Poisson brackets for the OPs (see Theorem~\ref{P:OP});
they were derived from the combinations of \cite{NenSim}, to which they are equivalent (see Proposition~\ref{P:cWP}),
but we prove them directly by induction. It was recently brought to our attention that results related to those reported here
were obtained in the early 1980s by P. Kulish, \cite{Kul}, using $r$-matrix methods. Together
with \cite{KakMug1}, \cite{KakMug2}, this is one of the first papers which investigate the Poisson bracket used in our work.

Before proceeding to our results, let us give some background on
orthogonal polynomials on the unit circle. For more details, we
direct the reader to \cite{Simon1} and \cite{Simon2}. The recurrence
relation \eqref{E:PhiRec} implies that
$\|\Phi_n\|_{L^2(d\mu)}=\prod_{j=0}^{n-1}\rho_j$, where
$\rho_j=\sqrt{1-|\al_j|^2}$. In particular, this means that
$\phi_n=\Phi_n/\prod_{j=0}^{n-1}\rho_j$ and the recurrence relations
for the orthonormal polynomials read
\begin{equation}\label{E:phiRec}
\rho_n\phi_{n+1}(z)=z\phi_n(z)-\bar\al_n\phi_n^*(z),
\end{equation}
\begin{equation}\label{E:phi*Rec}
\rho_n\phi_{n+1}^*(z)=\phi_n^*(z)-\al_nz\phi_n(z).
\end{equation}
Writing \eqref{E:PhiRec} and \eqref{E:Phi*Rec} in matrix form leads to
\begin{equation}\label{E:OPWall}
\begin{bmatrix}
\Phi_{n+1}(z)\\ \Phi_{n+1}^*(z)
\end{bmatrix}=
\begin{bmatrix}
z &  -\bar\al_n\\
-\al_nz   &  1
\end{bmatrix}
\cdot\begin{bmatrix}
\Phi_{n}(z)\\ \Phi_{n}^*(z)
      \end{bmatrix}
=\begin{bmatrix}
 zB_n^*(z) & -A_n^*(z) \\
 -zA_n(z)  & B_n(z)
 \end{bmatrix}
\cdot \begin{bmatrix}
\Phi_{0}(z)\\ \Phi_{0}^*(z)
      \end{bmatrix}.
\end{equation}
The $A_n$ and $B_n$ are monic polynomials of degree $n$, and they are known as Wall polynomials. For the proof
of the second identity in \eqref{E:OPWall}, see for example Theorem~3.2.10 of \cite{Simon1}.
Furthermore, if $\lambda$ is a complex unimodular number, and $\{\Phi_n^\lambda\}_{n\geq0}$ are the monic orthogonal polynomials associated to the (rotated) Verblunsky coefficients $\{\lambda\al_n\}_{n\geq0}$,
then the recurrence relations for these polynomials become
\begin{equation}\label{E:LRecRel}
\begin{aligned}
\begin{bmatrix}
\Phi_{n+1}^\lambda(z)\\ \bar\lambda(\Phi_{n+1}^\lambda)^*(z)
\end{bmatrix}
&=\begin{bmatrix}
 zB_n^*(z) & -A_n^*(z) \\
 -zA_n(z)  & B_n(z)
 \end{bmatrix}
\cdot \begin{bmatrix}
\Phi_{0}^\lambda(z)\\ \bar\lambda(\Phi_{0}^\lambda)^*(z)
      \end{bmatrix}\\
&=\begin{bmatrix}
 zB_n^*(z) & -A_n^*(z) \\
 -zA_n(z)  & B_n(z)
 \end{bmatrix}
\cdot \begin{bmatrix}
1\\ \bar\lambda
      \end{bmatrix}.
\end{aligned}
\end{equation}
The case $\lambda=-1$ is special, and $\Psi_n\equiv\Phi_n^{\lambda=-1}$ are called
second kind polynomials. Putting all the relations above together, one can easily obtain the Pinter-Nevai formulae, relating the monic orthogonal and second kind polynomials to the Wall polynomials:
\begin{eqnarray}
\Phi_n(z)&=&zB_{n-1}^*(z)-A_{n-1}^*(z)\label{E:PN1}\\
\Psi_n(z)&=&zB_{n-1}^*(z)+A_{n-1}^*(z).\label{E:PN11}
\end{eqnarray}
or, equivalently,
\begin{eqnarray}\label{E:PN2}
A_n(z)&=&\tfrac{1}{2z}\bigl(\Psi_{n+1}^*(z)-\Phi_{n+1}^*(z)\bigr)\\
B_n(z)&=&\tfrac12\bigl(\Psi_{n+1}^*(z)+\Phi_{n+1}^*(z)\bigr).\label{E:PN21}
\end{eqnarray}
The following relations will also be useful farther in our work:
\begin{equation}\label{E:PhiL}
\Phi_n^\lambda=\tfrac{1+\bar\lambda}{2}\cdot\Phi_n+\tfrac{1-\bar\lambda}{2}\cdot\Psi_n,\quad
\Psi_n^\lambda=\tfrac{1-\bar\lambda}{2}\cdot\Phi_n+\tfrac{1+\bar\lambda}{2}\cdot\Psi_n.
\end{equation}

\textit{Notation.} In order to keep our notation straight, we specify it here.
For a sequence of Verblunsky coefficients $\{\al_n\}_{n\geq0}$, we denote
\begin{eqnarray*}
\Phi_n(z)&\equiv&\Phi_n(z;\al_0,\ldots,\al_{n-1}),\\
\Psi_n(z)&\equiv&\Phi_n(z;-\al_0,\ldots,-\al_{n-1}),\\
\Phi_n^\lambda(z)&\equiv&\Phi_n(z;\lambda\al_0,\ldots,\lambda\al_{n-1})\\
\Psi_n^\lambda(z)&\equiv&\Phi_n(z;-\lambda\al_0,\ldots,-\lambda\al_{n-1}).
\end{eqnarray*}
Throughout this paper, we think of these polynomials, as well as of the Wall polynomials $A_n$
and $B_n$, as functions of the Verblunsky coefficients, with $z\in\IC$ a parameter. A simple calculation
using the recurrence relations shows that, for any $n\geq 0$:
$$
\|\Phi_n\|=\|\Phi_n^*\|=\|\Psi_n\|=\|\Psi_n^*\|=\prod_{j=0}^{n-1} \rho_j,
$$
where all the norms are taken in $L^2(d\mu)$, and for $n=0$ we use the usual convention that
the empty product is identically equal to 1. So the orthonormal and normalized second kind polynomials
are given by
$$
\phi_n=\Bigl(\prod_{j=0}^{n-1} \rho_j^{-1}\Bigr)\cdot \Phi_n\quad\text{and}\quad
\psi_n=\Bigl(\prod_{j=0}^{n-1} \rho_j^{-1}\Bigr)\cdot \Psi_n\,.
$$

\begin{remark}
As observed in \cite{CanSim}, one may use the detailed knowledge of Poisson brackets of
the first and second kind orthogonal polynomials to compute other, very interesting, Poisson brackets.
For example, it was shown in the work of Gekhtman and Nenciu, \cite{GekNen}, that in the finite case,
the whole Poisson structure can be encoded in the Poisson bracket of the associated Carath\'eodory function
at two arbitrary points $z\neq w\in\IC\setminus S^1$:
\begin{equation}\label{E:Carat}
\{F(z),F(w)\}=i\bigl(F(z)-F(w)\bigr)\bigl(F(z)F(w)-1\bigr)-i\frac{z+w}{z-w}\bigl(F(z)-F(w)\bigr)^2\,.
\end{equation}
The first proof of this formula (see \cite{GekNen} for some of its very interesting consequences) was
done in \cite{GekNen} using the r-matrix formulation of the Ablowitz-Ladik bracket (see also \cite{KilNen} 
and \cite{Li} for more details on this alternate way of defining the Poisson bracket \eqref{E:PbDefn}). 

Cantero and Simon show in \cite{CanSim} that one can more easily deduce this formula directly from the brackets
of the orthogonal polynomials. Indeed, in the finite case they use an analogue for measures with finite support of
Theorem~3.2.4 of \cite{Simon1} to connect the Carath\'eodory function with the orthogonal
polynomials of the first and second kind. Given the importance of formula~\eqref{E:Carat} to the study of
the finite Ablowitz-Ladik equation (see \cite{GekNen}), one would like to find the analogue
of \eqref{E:Carat} in the periodic case. This is particularly interesting right now in view of the results
and the open questions
of Li and Nenciu \cite{LiNen} for the periodic Ablowitz-Ladik equation. In order to go from the formulae
in our Theorem~\ref{P:OP} to a bracket for the Carath\'eodory function, one can in principle
use formulae (11.3.5)--(11.3.8) from \cite{Simon2}, which show that $F$ is the solution
of a quadratic equation whose coefficients are combinations of the first and second kind
polynomials. This leads to a straightforward algorithm for computing $\{F(z),F(w)\}$, but what
is neither clear nor easy is to ``close'' the formula, i.e. to express the result only
in terms of $F(z)$ and $F(w)$.  
\end{remark}

\section{Many Poisson brackets}

The main statement of this note is the following:
\begin{theorem}\label{P:OP}
Let $n\geq0$ and $d\mu$ be a probability measure on the circle, and such that either the support of $d\mu$ is 
finite and contains at least $n+1$ points, or the Verblunsky coefficients of $d\mu$ are periodic with 
minimal period $p\geq n+1$.
Then the Poisson brackets\footnote{Using either definition \eqref{E:PbDefn} or \eqref{E:PbDefnPeriodic}, as appropriate.} of the monic orthogonal and second kind polynomials associated to this measure are given, for any $z\neq w\in\IC$, by
\begin{equation}\label{E:PhiPhi}
\{\Phi_n(z),\Phi_n(w)\}=\{\Psi_n(z),\Psi_n(w)\}=0,
\end{equation}
\begin{equation}\label{E:PhiPhi*}
\{\Phi_n(z),\Phi^*_n(w)\}=\tfrac{iw}{z-w}\bigl(\Phi_n(z)\Phi_n^*(w)-\Phi_n(w)\Phi_n^*(z)\bigr),
\end{equation}
\begin{equation}\label{E:PsiPsi*}
\{\Psi_n(z),\Psi_n^*(w)\}=\tfrac{iw}{z-w}\bigl(\Psi_n(z)\Psi_n^*(w)-\Psi_n(w)\Psi_n^*(z)\bigr),
\end{equation}
\begin{equation}\label{E:PhiPsi}
\begin{aligned}
\{\Phi_n(z),\Psi_n(w)\}
&=\tfrac{-iw}{z-w}\bigl(\Phi_n(z)\Psi_n(w)-\Phi_n(w)\Psi_n(z)\bigr)\\
&\quad -\tfrac{i}{2}\bigl(\Phi_n(z)-\Psi_n(z)\bigr)\bigl(\Phi_n^*(w)+\Psi_n^*(w)\bigr),
\end{aligned}
\end{equation}
\begin{equation}\label{E:PhiPsi*}
\begin{aligned}
\{\Phi_n(z),\Psi_n^*(w)\}
&=\tfrac{-iw}{z-w}\bigl(\Psi_n(z)\Phi^*_n(w)-\Psi_n(w)\Phi^*_n(z)\bigr)\\
&\quad +\tfrac{i}{2}\bigl(\Phi_n(z)-\Psi_n(z)\bigr)\bigl(\Phi_n^*(w)-\Psi_n^*(w)\bigr).
\end{aligned}
\end{equation}
\end{theorem}

\begin{remark}\label{R:OP}
Note that, by using the invariance of the Poisson bracket under complex conjugation
\begin{equation}\label{E:barPB}
\{\bar f,\bar g\}=\overline{\{f,g\}}
\end{equation}
and the formula for reversed polynomials, one can compute the remaining Poisson brackets from the ones
in Theorem~\ref{P:OP}. We obtain
\begin{equation}\label{E:P*P*}
\{\Phi^*_n(z),\Phi^*_n(w)\}=\{\Psi^*_n(z),\Psi^*_n(w)\}=0,
\end{equation}
\begin{equation}\label{E:Phi*Psi}
\begin{aligned}
\{\Phi_n^*(z),\Psi_n(w)\}
&=\tfrac{-iz}{z-w}\bigl(\Psi_n^*(z)\Phi_n(w)-\Psi_n^*(w)\Phi_n(z)\bigr)\\
&\quad -\tfrac{i}{2}\bigl(\Phi_n^*(z)-\Psi_n^*(z)\bigr)\bigl(\Phi_n(w)-\Psi_n(w)\bigr),
\end{aligned}
\end{equation}
and
\begin{equation}\label{E:Phi*Psi*}
\begin{aligned}
\{\Phi_n^*(z),\Psi_n^*(w)\}
&=\tfrac{-iz}{z-w}\bigl(\Phi_n^*(z)\Psi_n^*(w)-\Phi_n^*(w)\Psi_n^*(z)\bigr)\\
&\quad +\tfrac{i}{2}\bigl(\Phi_n^*(z)-\Psi_n^*(z)\bigr)\bigl(\Phi_n(w)+\Psi_n(w)\bigr).
\end{aligned}
\end{equation}

Also note that, since the left and right-hand sides of the all the formulae given above, \eqref{E:PhiPhi}--\eqref{E:Phi*Psi*},
are continuous in $z,w\in\IC$, we can find the Poisson brackets for $z=w$ in the usual fashion. For example, formula \eqref{E:PhiPsi*}
becomes, for $z=w$,
\begin{equation}\label{E:PhiPsi*2}
\begin{aligned}
\{\Phi_n(z),\Psi_n^*(z)\}
&=-iz\bigl(\Psi_n^\prime(z)\Phi^*_n(z)-\Psi_n(z)(\Phi^*_n)^\prime(z)\bigr)\\
&\quad +\tfrac{i}{2}\bigl(\Phi_n(z)-\Psi_n(z)\bigr)\bigl(\Phi_n^*(z)-\Psi_n^*(z)\bigr),
\end{aligned}
\end{equation}
where $^\prime$ denotes the usual complex differentiation.
\end{remark}

\begin{remark}
Most of the Poisson brackets from Theorem~\ref{P:OP} and Remark~\ref{R:OP}
have already been computed by Cantero and Simon \cite{CanSim}. The Poisson bracket that they didn't compute is \eqref{E:PhiPsi*} (and its
reverse, \eqref{E:Phi*Psi}).
With this, we know all the brackets of orthogonal polynomials on the unit circle.
\end{remark}

\begin{proof}
First note that the polynomials of first and second kind and their reverses $\Phi_n,\Phi_n^*,\Psi_n,\Psi_n^*$
depend only on $\alpha_0,...,\alpha_{n-1}$, and so the hypothesis of the theorem imply that
the expressions of their Poisson brackets in the finite and periodic cases will be identical. We can thus assume
wlog that we are in the finite case.

For $n=0$, all the polynomials are identically equal to 1, and the relations hold trivially.
Assume all the identities hold for some $n\geq0$. We will prove \eqref{E:PhiPsi*} for $n+1$. All the
other relations follow similarly.

Assume that $|\text{supp}(d\mu)|>n+1$ (and hence $\al_0,\dots,\al_n\in \ID$). Using the recurrence relations, we immediately get
\begin{equation}\label{E:PP*Calc}
\begin{aligned}
\{\Phi_{n+1}(z),\Psi_{n+1}^*(w)\}
&=\{z\Phi_n(z)-\bar\al_n\Phi_n^*(z), \Psi_n^*(w)+\al_n w\Psi_n(w)\}\\
&\begin{aligned}
 = & z\{\Phi_n(z),\Psi_n^*(w)\}+\al_n zw\{\Phi_n(z),\Psi_n(w)\}\\
   &  -\bar\al_n\{\Phi_n^*(z),\Psi_n^*(w)\}-|\al_n|^2w\{\Phi_n^*(z),\Psi_n(w)\}\\
   & -\{\bar\al_n,\al_n\} w\Phi_n^*(z)\Psi_n(w)\\
 \end{aligned}
\end{aligned}
\end{equation}
The last identity follows since the $n^\text{th}$ polynomials depend only on $\al_0,\dots,\al_{n-1}$,
and hence their Poisson brackets with $\al_n$ or $\bar\al_n$ are identically 0. Using the induction
hypothesis, we substitute the formulae for the brackets on the right hand side and we get
\begin{equation*}
\{\Phi_{n+1}(z),\Psi_{n+1}^*(w)\}=\tfrac{-iw}{z-w}\cdot T_1 +\tfrac{i}{2}\cdot T_2 -i\rho_n^2w\Phi_n^*(z)\Psi_n(w),
\end{equation*}
where
\begin{equation*}
\begin{aligned}
T_1 &=\bigl[\Psi_{n+1}(z)\Phi_{n+1}^*(w)-\Psi_{n+1}(w)\Phi_{n+1}^*(z)\bigr]\\
    &\quad +(w-z)\bigl(|\al_n|^2\Phi_n(w)\Psi_n^*(z)+\Psi_n(w)\Phi_n^*(z)\bigr)
\end{aligned}
\end{equation*}
and
\begin{equation*}
\begin{aligned}
T_2 &=z\bigl(\Phi_n(z)-\Psi_n(z)\bigr)\bigl(\Phi_n^*(w)-\Psi_n^*(w)\bigr)
    -\bar\al_n \bigl(\Phi_n^*(z)+\Psi_n^*(z)\bigr)\bigl(\Phi_n^*(w)-\Psi_n^*(w)\bigr)\\
    &\quad -\al_nzw\bigl(\Phi_n(z)-\Psi_n(z)\bigr)\bigl(\Phi_n(w)+\Psi_n(w)\bigr)
    +|\al_n|^2w\bigl(\Phi_n^*(z)-\Psi_n^*(z)\bigr)\bigl(\Phi_n(w)-\Psi_n(w)\bigr)\\
    &=\bigl[\Phi_{n+1}(z)-\Psi_{n+1}(z)\bigr]\bigl[\Phi_{n+1}^*(w)-\Psi_{n+1}^*(w)\bigr]
    -2|\al_n|^2w\bigl(\Phi_n^*(z)\Psi_n(w)+\Phi_n(w)\Psi_n^*(z)\bigr).
\end{aligned}
\end{equation*}
If we add together all the extra terms, we get
\begin{align*}
&iw\bigl(|\al_n|^2\Phi_n(w)\Psi_n^*(z)+\Psi_n(w)\Phi_n^*(z)\bigr)
-i|\al_n|^2w\bigl(\Phi_n^*(z)\Psi_n(w)+\Phi_n(w)\Psi_n^*(z)\bigr)\\
&-i(1-|\al_n|^2)w\Phi_n^*(z)\Psi_n(w)\equiv0,
\end{align*}
which proves \eqref{E:PhiPsi*} for $n+1$ and $|\text{supp}(d\mu)|>n+1$.

The proofs of the other relations \eqref{E:PhiPhi}--\eqref{E:PhiPsi}
run along the same lines.
\end{proof}

\begin{remark}
Note that, if in the second step of the induction performed above, we assume that 
$|\text{supp}(d\mu)|=n+1$, then $\al_n\in S^1$ and
it is a Casimir for the corresponding Poisson bracket. In particular $\{\bar\al_n,\al_n\}=0$ in the 
last term of \eqref{E:PP*Calc}. Keeping this in mind,
and using the fact that $|\al_n|=1$, one can still run through the calculation above and get the result.
By induction, \eqref{E:PhiPsi*} holds
for all $n\geq0$, as long as $|\text{supp}(d\mu)|\geq n$. That shows that the statements
of the theorem, \eqref{E:PhiPhi}--\eqref{E:PhiPsi*}, hold even if $|\text{supp}(d\mu)|=n$.
\end{remark}

We use the results of Theorem~\ref{P:OP} to find the Poisson brackets
of the orthonormal and normalized second kind orthogonal polynomials $\phi_n$, $\psi_n$, and their
reverses $\phi_n^*$ and $\psi_n^*$.

\begin{lemma}\label{L:rhoPb}
Let $n\geq 0$, and define $R_n=\prod_{k=0}^{n-1} \rho_k^{-1}$ if $n\geq 1$, $R_0=1$, to be the inverse
of the norm in $L^2(d\mu)$ of the degree $n$ polynomials. Then
\begin{equation}\label{E:rhoPhi}
-\{R_n,\Phi_n(z)\}=\{R_n,\Psi_n(z)\}=\tfrac{i}{4R_n}\bigl(\Phi_n(z)-\Psi_n(z)\bigr)
=\tfrac{i}{4}\bigl(\phi_n(z)-\psi_n(z)\bigr),
\end{equation}
and
\begin{equation}\label{E:rhoPhi*}
\{R_n,\Phi_n^*(z)\}=-\{R_n, \Psi_n^*(z)\}=\tfrac{i}{4R_n}\bigl(\Phi_n^*(z)-\Psi_n^*(z)\bigr)
=\tfrac{i}{4}\bigl(\phi_n^*(z)-\psi_n^*(w)\bigr).
\end{equation}
\end{lemma}

\begin{proof}
For $0\geq j\geq n-1$ one has $\{R_n,\al_j\}=\tfrac{i}{2}R_n\al_j$,
and hence the evolution of the $\al_j$'s under the flow generated by the Hamiltonian $R_n$ ($n>j$) is given by
$$
\al_j(t)=\al_je^{i\frac{t}{2R_n}}.
$$
If we denote $e^{i\frac{t}{2R_n}}=\lambda(t)$, with $|\lambda(t)|=1$, we obtain using \eqref{E:PhiL} that:
\begin{align*}
\{R_n,\Phi_n(z)\}
&=\left.\frac{d}{dt}\right|_{t=0} \Phi_n\bigl(z; \{\lambda(t)\al_j\}_{j=0}^{n-1}\bigr)
=\left.\frac{d}{dt}\right|_{t=0} \Bigl(\tfrac{1+\overline{\lambda(t)}}{2}\Phi_n(z) + \tfrac{1-\overline{\lambda(t)}}{2}\Psi_n(z) \Bigr)\\
&=-\frac{i}{4R_n}\bigl(\Phi_n(z)-\Psi_n(z)\bigr)
\end{align*}
The other three Poisson brackets are computed in a similar fashion, using \eqref{E:PhiL}
and the analogous statements for the reversed polynomials.
\end{proof}

This lemma allows us to compute the Poisson brackets of the normalized OP's:
\begin{prop}\label{P:NOP}
If $n\geq0$ and $|\text{supp}(d\mu)|\geq n+1$ or the Verblunsky coefficients of
$d\mu$ are periodic with period $p\geq n+1$, then one has the following Poisson brackets for $z\neq w\in\IC$:
\begin{equation}\label{E:phiphi}
\{\phi_n(z),\phi_n(w)\}=-\{\psi_n(z),\psi_n(w)\}=\tfrac{i}{4}\bigl(\phi_n(z)\psi_n(w)-
\phi_n(w)\psi_n(z)\bigr),
\end{equation}
\begin{equation}\label{E:phiphi*}
\begin{aligned}
\{\phi_n(z),\phi_n^*(w)\}
&=\tfrac{i}{4}\bigl[\phi_n(z)(\phi_n^*(w)-\psi_n^*(w))+(\phi_n(z)-\psi_n(z))\phi_n^*(w)\bigr]\\
&\quad +\tfrac{iw}{z-w}\bigl(\phi_n(z)\phi_n^*(w)-\phi_n(w)\phi_n^*(z)\bigr),
\end{aligned}
\end{equation}
\begin{equation}\label{E:psipsi*}
\begin{aligned}
\{\psi_n(z),\psi_n^*(w)\}
&=\tfrac{i}{4}\bigl[\psi_n(z)(\psi_n^*(w)-\phi_n^*(w))+(\psi_n(z)-\phi_n(z))\psi_n^*(w)\bigr]\\
&\quad +\tfrac{iw}{z-w}\bigl(\psi_n(z)\psi_n^*(w)-\psi_n(w)\psi_n^*(z)\bigr),
\end{aligned}
\end{equation}
\begin{equation}\label{E:phipsi}
\begin{aligned}
\{\phi_n(z),\psi_n(w)\}
&=\tfrac{-iw}{z-w}\bigl(\phi_n(z)\psi_n(w)-\phi_n(w)\psi_n(z)\bigr)\\
&\quad-\tfrac{i}{2}\bigl(\phi_n(z)-\psi_n(z)\bigr)\bigl(\phi_n^*(w)+\psi_n^*(w)\bigr)\\
&\quad+\tfrac{i}{4}\bigl(\phi_n(z)\phi_n(w)-\psi_n(z)\psi_n(w)\bigr),
\end{aligned}
\end{equation}
\begin{equation}\label{E:phipsi*}
\begin{aligned}
\{\phi_n(z),\psi_n^*(w)\}
&=\tfrac{i}{4}\bigl[(\phi_n(z)-\psi_n(z))\phi_n^*(w)-\psi_n(z)(\phi_n^*(w)-\psi_n^*(w))\bigr]\\
&\quad -\tfrac{iw}{z-w}\bigl(\psi_n(z)\phi_n^*(w)-\psi_n(z)\phi_n^*(z)\bigr).
\end{aligned}
\end{equation}
\end{prop}

\begin{remark}
Same as in the case of the monic polynomials (Theorem~\ref{P:OP}), one can obtain the Poisson brackets
for $z=w$ by taking limits as $z\to w$ in the formulae for $z\neq w$.

Note that the condition on the size of the support of the measure $d\mu$ is more restrictive in this proposition than
in Theorem~\ref{P:OP}. Indeed, if $|\text{supp}(d\mu)|=n$, then one can define the monic polynomial $\Phi_n$ by
the usual orthogonalization procedure applied to $z^n$, with respect to $\Phi_0,\dots,\Phi_{n-1}$. But since
$\Phi_0,\dots,\Phi_{n-1}$ form a basis in $L^2(d\mu)$, this
leads to a monic polynomial of degree $n$ which is identically zero in $L^2(d\mu)$. In particular, one cannot normalize
this polynomial, and hence the only orthonormal polynomials that can be correctly defined are $\phi_0,\dots,\phi_{n-1}$.
\end{remark}

\begin{proof}
The results of this proposition follow by simple calculations from Lemma~\ref{L:rhoPb} and
Theorem~\ref{P:OP}, since $R_n$ is the normalization constant for both the regular and second kind polynomials, and their reverses. For example, the first Poisson bracket \eqref{E:phiphi} can be computed as follows:
\begin{align*}
\{\phi_n(z),\phi_n(w)\}
=\{R_n,\Psi_n(w)\}R_n\Phi_n(z)+\{\Phi_n(z),R_n\}R_n\Psi_n(w)
\end{align*}
since $\{R_n,R_n\}$ and $\{\Phi_n(z),\Phi_n(w)\}$ are identically zero. Therefore,
using \eqref{E:rhoPhi} and canceling terms, we obtain \eqref{E:phiphi}. The other formulae
can be obtained similarly.
\end{proof}

Finally, we have the following:
\begin{prop}\label{P:WP}
The Poisson brackets of the Wall polynomials are given, for $z\neq w\in\IC$, by
\begin{equation}\label{E:AA}
\{A_n(z),A_n(w)\}=\{B_n(z),B_n(w)\}=0,
\end{equation}
\begin{equation}\label{E:AA*}
\{A_n(z),A_n^*(w)\}=iA_n(z)A_n^*(w)-\tfrac{i}{z-w}\bigl(zB_n^*(z)B_n(w)-wB_n^*(w)B_n(z)\bigr),
\end{equation}
\begin{equation}\label{E:BB*}
\{B_n(z),B_n^*(w)\}=\tfrac{iz}{z-w}\bigl(A_n(z)A_n^*(w)-A_n(w)A_n^*(z)\bigr),
\end{equation}
\begin{equation}\label{E:AB}
\{A_n(z),B_n(w)\}=\tfrac{-iw}{z-w}\bigl(A_n(z)B_n(w)-A_n(w)B_n(z)\bigr),
\end{equation}
and
\begin{equation}\label{E:AB*}
\{A_n(z),B_n^*(w)\}=\tfrac{iz}{z-w}\bigl(A_n(z)B_n^*(w)-A_n(w)B_n^*(z)\bigr).
\end{equation}
\end{prop}

\begin{remark}
All the other brackets for Wall polynomials can be obtained from the ones above by taking
the reversed polynomials in $z$ and $w$, and using \eqref{E:barPB}. The Poisson brackets for $z=w$ can
be obtained from the ones above by taking limits on both sides of the identities (see also Remark~2.1).
The six Poisson brackets \eqref{E:AA}--\eqref{E:AB*} for Wall polynomials are equivalent
to the Poisson brackets \eqref{E:PhiPhi}--\eqref{E:PhiPsi*}.
\end{remark}

\begin{proof}
The six Poisson brackets in this proposition are equivalent to the
six Poisson brackets for monic and second kind orthogonal
polynomials from Theorem~\ref{P:OP}. This equivalence is merely
a reflection of the Pinter-Nevai formulae \eqref{E:PN1}--\eqref{E:PN21}, relating orthogonal and
Wall polynomials.
The rest is just simple algebra.
\end{proof}

As explained in the Introduction, the first, to the best of our knowledge,
computations of Poisson brackets for polynomials associated to OPUC were done
in the work of Nenciu and Simon~\cite{NenSim} (see, also, \cite{Simon2} and \cite{Nen2}). The authors were interested in
understanding the Ablowitz-Ladik system in the periodic setting, and hence wanted
to find Poisson brackets for the discriminant associated to this problem. In terms of Wall polynomials, they
have to show that
\begin{equation*}
\{B_n(z)+zB_n^*(z),B_n(w)+wB_n^*(w)\}=0.
\end{equation*}
The main observation is that, since the $\al$'s are separated in the definition of
the Poisson bracket, one can use induction in $n$ to prove statements such as the one above.
Using the recurrence relations for Wall polynomials immediately shows that one needs to
compute certain anti-symmetric combinations of Poisson brackets involving both Wall polynomials ($A$ and $B$), and their
reverses. What they show is the following: For $n\geq0$, set
\begin{eqnarray*}
  F_n(z,w)=-i\{A_n^*(z),B_n(w)\},  &\qquad&
  R_n(z,w)=i\{B_n^*(z),B_n(w)\}, \\
  S_n(z,w)=i\{A_n(z),A_n^*(w)\},   &\qquad&
  X_n(z,w)=i\{B_n^*(z),A_n^*(w)\}.
\end{eqnarray*}
and, for $q\in\ZZ$ and $z,w\in\IC\setminus\{0\}$, define
\begin{equation*}\label{defnQ}
Q_q(z,w)=zw\cdot\frac{z^{q-1}-w^{q-1}}{z-w}\,.
\end{equation*}
Condensing the statements of Nenciu and Simon (see Propositions 4.5 and 4.7 of \cite{Nen2}), one has
\begin{prop}\label{P:cWP}
The following statements hold for all $n\geq0$, $q\in\ZZ$, and
$z,w\neq0$:
\begin{eqnarray*}
  \emph{(}\al_n,\beta_n\emph{)} &\quad& \{A_n^*(z),A_n^*(w)\}=\{B_n(z),B_n(w)\}=0 \\
  \emph{(}\gamma_{n,q}\emph{)} &\quad& z^qF_n(z,w)-w^qF_n(w,z)=Q_q(z,w)[A_n^*(z)B_n(w)-A_n^*(w)B_n(z)]\\
  (r_{n,q}) &\quad& z^qR_n(z,w)-w^qR_n(w,z)=Q_q(z,w)\big[A_n(z)A_n^*(w)-A_n(w)A_n^*(z)\big]\\
  (s_{n,q}) &\quad& \begin{aligned}z^qS_n(z,w)-w^qS_n(w,z)&=-\big[z^qA_n(z)A_n^*(w)-w^qA_n(w)A_n^*(z)\big]\\
  &+\frac{z^q-w^q}{z-w}\big[zB_n^*(z)B_n(w)-wB_n^*(w)B_n(z)\big]\\
  \end{aligned}\\
  (x_{n,q}) &\quad& z^qX_n(z,w)-w^qX_n(w,z)=Q_q(z,w)\big[B_n^*(z)A_n^*(w)-B_n^*(w)A_n^*(z)\big]\\
\end{eqnarray*}
\end{prop}

\begin{proof}
The original proof of this proposition goes by induction in $n$, where the $n^\text{th}$ induction statement consists of the identities
($\al_n$)--($x_{n,q}$) for all $q\in\ZZ$. This makes for a fairly complicated proof.

In fact, the sequence of identities $(\gamma_{n,q})_{q\in\ZZ}$ is equivalent to formula \eqref{E:AB*}.
Indeed, take $(\gamma_{n,q})$ for $q=0$ and 1 and insert the definition of $F_n$. A simple calculation shows
that $Q_0(z,w)\equiv -1$ and $Q_1(z,w)\equiv 0$, and hence
\begin{equation*}
\{A_n^*(z),B_n(w)\}-\{A_n^*(w),B_n(z)\}=-i\bigl(A_n^*(z)B_n(w)-A_n^*(w)B_n(z)\bigr)
\end{equation*}
\begin{equation*}
z\{A_n^*(z),B_n(w)\}-w\{A_n^*(w),B_n(z)\}=0
\end{equation*}
Multiplying the first equation by $w$ and subtracting yields
\begin{equation}\label{E:A*B}
\{A_n^*(z),B_n(w)\}=\tfrac{iw}{z-w}\bigl(A_n^*(z)B_n(w)-A_n^*(w)B_n(z)\bigr),
\end{equation}
which is equivalent to \eqref{E:AB*}. Conversely, consider the combination
on the left-hand side of $(\gamma_{n,q})$ for some $q\in\ZZ$ and insert \eqref{E:A*B}.
Since the factor
$$A_n^*(z)B_n(w)-A_n^*(w)B_n(z)$$ is anti-symmetric in $z$ and $w$, the
right-hand side of $(\gamma_{n,q})$ emerges from
\begin{equation*}
z^q\frac{w}{z-w}-w^q\frac{z}{w-z}\cdot(-1)=zw\frac{z^{q-1}-w^{q-1}}{z-w}\equiv Q_q(z,w).
\end{equation*}
Equally elementary considerations show that ($r_{n,q}$) through ($x_{n,q}$) are equivalent to
the corresponding identities from Proposition~\ref{P:WP}.
\end{proof}


\end{document}